\magnification=1200
\input amssym.def
\input amssym.tex

\hsize=13.5truecm
\baselineskip=16truept
\baselineskip=16truept
\font\secbf=cmb10 scaled 1200
\font\eightrm=cmr8
\font\sixrm=cmr6

\font\eighti=cmmi8

\font\sixi=cmmi6
\skewchar\eighti='177 \skewchar\sixi='177

\font\eightsy=cmsy8
\font\sixsy=cmsy6
\skewchar\eightsy='60 \skewchar\sixsy='60

\font\eightit=cmti8

\font\eightbf=cmbx8
\font\sixbf=cmbx6

\let\sc=\tensc

\font\eightsc=cmcsc10 scaled 800
\font\secbf=cmb10 scaled 1200
\font\subsecfont=cmb10 scaled \magstephalf
\font\amb=cmmib10

\font\ambi=cmmib10 scaled 700

\newfam\mbfam 

\textfont\mbfam\amb \scriptfont\mbfam\ambi


\def\aa{\def\rm{\fam0\eightrm}%
  \textfont0=\eightrm \scriptfont0=\sixrm \scriptscriptfont0=\fiverm
  \textfont1=\eighti \scriptfont1=\sixi \scriptscriptfont1=\fivei
  \textfont2=\eightsy \scriptfont2=\sixsy \scriptscriptfont2=\fivesy
  \textfont3=\tenex \scriptfont3=\tenex \scriptscriptfont3=\tenex
  \def\sc{\eightsc}
  \def\it{\fam\itfam\eightit}%
  \textfont\itfam=\eightit
  \def\bf{\fam\bffam\eightbf}%
  \textfont\bffam=\eightbf \scriptfont\bffam=\sixbf
   \scriptscriptfont\bffam=\fivebf
  \normalbaselineskip=9.7pt
  \setbox\strutbox=\hbox{\vrule height7pt depth2.6pt width0pt}%
  \normalbaselines\rm}

\def\Proof{\vskip12pt\noindent{\bf Proof.} }

\def\Def#1{\vskip12pt\noindent{\bf Definition #1}}
\def\Remark#1{\vskip12pt\noindent{\bf Remark #1}}

\def\m@th{\mathsurround=0pt}

\def\cc#1{\hbox to .89\hsize{$\displaystyle\hfil{#1}\hfil$}\cr}
\def\lc#1{\hbox to .89\hsize{$\displaystyle{#1}\hfill$}\cr}
\def\rc#1{\hbox to .89\hsize{$\displaystyle\hfill{#1}$}\cr}

\def\eqal#1{\null\,\vcenter{\openup\jot\m@th
  \ialign{\strut\hfil$\displaystyle{##}$&&$\displaystyle{{}##}$\hfil
      \crcr#1\crcr}}\,}

\def\section#1{\vskip 22pt plus6pt minus2pt\penalty-400
        {{\secbf
        \noindent#1\rightskip=0pt plus 1fill\par}}
        \par\vskip 12pt plus5pt minus 2pt
        \penalty 1000}

\def\subsection#1{\vskip 20pt plus6pt minus2pt\penalty-400
        {{\subsecfont
        \noindent#1\rightskip=0pt plus 1fill\par}}
        \par\vskip 8pt plus5pt minus 2pt
        \penalty 1000}

\def\subsubsection#1{\vskip 18pt plus6pt minus2pt\penalty-400
        {{\subsecfont
        \noindent#1}}
        \par\vskip 7pt plus5pt minus 2pt
        \penalty 1000}

\def\center#1{{\begingroup \leftskip=0pt plus 1fil\rightskip=\leftskip
\parfillskip=0pt \spaceskip=.3333em \xspaceskip=.5em \pretolerance 9999
\tolerance 9999 \parindent 0pt \hyphenpenalty 9999 \exhyphenpenalty 9999
\par #1\par\endgroup}}

\def\\{\hfill\break}

\def\kwadrat{\null\ \hfill\null\ \hfill$\square$}
\def\mida#1{{{\null\kern-4.2pt\left\bracevert\vbox to 6pt{}\!\hbox{$#1$}\!\right\bracevert\!\!}}}
\def\midy#1{{{\null\kern-4.2pt\left\bracevert\!\!\hbox{$\scriptstyle{#1}$}\!\!\right\bracevert\!\!}}}

\def\diagint{{\raise1.5pt\hbox{$\scriptscriptstyle\diagup$}\hskip-8.7pt\intop}}

\def\divv{{\rm div}\,}

\def\today{${\scriptscriptstyle\number\day-\number\month-\number\year}$}
\footline={{\hfil\rm\the\pageno\hfil${\scriptscriptstyle\rm\jobname}$\ \ \today}}
\def\D{{\Bbb D}}

\def\R{{\Bbb R}}

\def\meas{{\rm meas\,}}
\def\esssup{\mathop{\rm esssup\,}\limits}

\center{\secbf The H\"older continuity of the swirl for the Navier-Stokes 
motions}
\vskip1.5cm

\centerline{\bf Wojciech M. Zaj\c aczkowski}

\vskip1cm
\noindent
Institute of Mathematics, Polish Academy of Sciences,\\
\'Sniadeckich 8, 00-956 Warsaw, Poland\\
E-mail:wz@impan.pl;\\
Institute of Mathematics and Cryptology, Cybernetics Faculty,\\
Military University of Technology, Kaliskiego 2,\\
00-908 Warsaw, Poland
\vskip0.8cm


\vskip1.5cm

\noindent
{\bf Abstract.} 
We examine the axially symmetric solutions to the Navier-Stokes equations in 
a periodic cylinder $\Omega\subset\R^3$. By swirl we denote the expression 
$u=rv_\varphi$, where $v_\varphi$ is the angular component of velocity and 
$r$, $\varphi$, $z$ are the cylindrical coordinates. The swirl satisfies 
a linear parabolic equation with coefficients equal to velocity. Assuming that 
$v\in L_q(0,T;L_p(\Omega))$, ${3\over p}+{2\over q}=1-{3\over2}\varkappa$, 
$\varkappa>0$, we show that $u\in C^{3\varkappa/2,3\varkappa/4}(\Omega^T)$, 
$T>0$. We apply the DeGiorgi method developed in Chapter 2 of the book: "Linear 
and quasilinear equations of parabolic type" written by Ladyzhenskaya, 
Solonnikov and Uraltseva.
\vfil\eject

\section{1. Introduction}

In this paper we examine some H\"older properties of axially symmetric 
solutions to the problem
$$\eqal{
&v_t+v\cdot\nabla v-\nu\Delta v+\nabla p=0\quad &{\rm in}\ \ 
\Omega^T=\Omega\times(0,T),\cr
&\divv v=0\quad &{\rm in}\ \ \Omega^T,\cr
&v\cdot\bar n=0,\ \ \bar n\cdot\D(v)\cdot\bar\tau_\alpha=0,\ \ \alpha=1,2,\quad
&{\rm on}\ \ S_1^T=S_1\times(0,T),\cr
&{\rm periodic\ boundary\ conditions}\quad &{\rm on}\ \ S_2^T,\cr
&v|_{t=0}=v_0\quad &{\rm in}\ \ \Omega,\cr}
\leqno(1.1)
$$
where $\Omega\subset\R^3$ is a cylinder with the boundary $S=S_1\cup S_2$, 
$S_1$ is parallel to the axis of the cylinder and $S_2$ is perpendicular 
to it.

\noindent
By $v=v(x,t)=(v_1(x,t),v_2(x,t),v_3(x,t))\in\R^3$ we denote velocity of the 
fluid, $p=p(x,t)\in\R$ the pressure, $x=(x_1,x_2,x_3)$ are the Cartesian 
coordinates such that $x_3$-axis is the axis of the cylinder $\Omega$. 
By $\bar n$ we denote the unit outward normal vector to $S_1$, 
$\bar\tau_\alpha$, $\alpha=1,2$, is the tangent vector to $S_1$. Moreover, 
the dot denotes the scalar product in $\R^3$, $\D(v)=\nabla v+\nabla v^T$ 
and $\nu>0$ is the constant viscosity coefficient.

\noindent
Since we are interested to examine axially symmetric solutions to (1.1) we 
introduce the cylindrical coordinates $r$, $\varphi$, $z$ by the relations
$$
x_1=r\cos\varphi,\quad x_2=r\sin\varphi,\quad x_3=z.
$$
Then $\Omega$ must be an axially symmetric cylinder.

\noindent
To be more precise we fix positive numbers $R$ and $a$ such that
$$
\Omega=\{x\in\R^3:\ r<R\ {\rm and}\ |z|<a\}
$$
and
$$\eqal{
&S_1=\{x\in\R^3:\ r=R,\ |z|\le a\},\cr
&S_2=\{x\in\R^3:\ r<R\ z\in\{-a,a\}\}.\cr}
$$
Let us introduce the vectors
$$
\bar e_r=(\cos\varphi,\sin\varphi,0),\quad
\bar e_\varphi=(-\sin\varphi,\cos\varphi,0),\quad \bar e_z=(0,0,1).
$$
Then the cylindrical components of velocity are described by the relations
$$
v_r=v\cdot\bar e_r,\quad v_\varphi=v\cdot\bar e_\varphi,\quad
v_z=v\cdot\bar e_z.
\leqno(1.2)
$$

\Def{1.1.}
The axially symmetric solutions to problem (1.1) are such that
$$
v_{r,\varphi}=v_{\varphi,\varphi}=v_{z,\varphi}=p_{,\varphi}=0.
\leqno(1.3)
$$

\Def{1.2.}
By the swirl we mean the quantity
$$
u=rv_\varphi.
\leqno(1.4)
$$
The aim of this paper is to show boundedness and the H\"older continuity of 
$u$. Especially, we are interested to estimate the H\"older exponent. For this 
purpose we have to repeat the considerations from [LSU, Ch. 2, Sect. 5--8].

From [LL, Ch. 2] and [Z2, Ch. 2] it follows the following problem for 
$v_\varphi$
$$\eqal{
&v_{\varphi,t}+v\cdot\nabla v_\varphi+{v_r\over r}v_\varphi-\nu\Delta v_\varphi
+\nu{v_\varphi\over r^2}=0\quad &{\rm in}\ \ \Omega^T,\cr
&v_{\varphi,r}={1\over R}v_\varphi\quad &{\rm on}\ \ S_1^T,\cr
&{\rm periodic\ boundary\ condition}\quad &{\rm on}\ \ S_2^T,\cr
&v_\varphi|_{t=0}=v_\varphi(0)\quad &{\rm in}\ \ \Omega,\cr}
\leqno(1.5)
$$
where $\Delta v_\varphi={1\over r}(rv_{\varphi,r})_{,r}+v_{\varphi,zz}$, 
$v\cdot\nabla=v_r\partial_r+v_z\partial_z$.

\noindent
In view of (1.4) problem (1.5) implies the following problem for $u$,
$$\eqal{
&u_t+v\cdot\nabla u-\nu\Delta u+\nu{u_{,r}\over r}=0\quad &{\rm in}\ \ 
\Omega^T,\cr
&u_{,r}={2\over R}u\quad &{\rm on}\ \ S_1^T,\cr
&{\rm periodic\ boundary\ condition}\quad &{\rm on}\ \ S_2^T,\cr
&u|_{t=0}=u_0\quad &{\rm in}\ \ \Omega.\cr}
\leqno(1.6)
$$
Let $V_2^0(\Omega^T)$ be a space of functions with the finite norm
$$
\|u\|_{L_\infty(0,T;L_2(\Omega))}+\|\nabla u\|_{L_2(\Omega^T)}.
$$
From [Z1, Lemma 2.1 and Remark 2.5] we have

\Remark{1.3.} 
For weak solutions to problem (1.1) (see [Z1, Z3]) we have
$$
\|v\|_{V_2^0(\Omega^T)}+\bigg\|{u\over r^2}\bigg\|_{L_2(\Omega^T)}\le 
c\|v_0\|_{L_2(\Omega)}\equiv d_1,
\leqno(1.7)
$$
where $d_1$ is some positive constant.

\proclaim Main Theorem. 
Assume that $v'\in L_{10}(\Omega^T)$, $v'=(v_r,v_z)$. 
Assume also that $u_0=rv_0\in L_\infty(\Omega)$. Assume that in 
a neighborhood $U$ of the axis of symmetry $u_0\in C^\alpha(\Omega)$, 
$\alpha={3\over2}\varkappa$, $\varkappa\in\big(0,{1\over3}\big]$. Then 
$u\in C^{\alpha,\alpha/2}(U\times(0,T))$ and its norm does not depend on $v'$.

In this paper we follow the idea of DeGiorgi developed in the book [LSU, 
Ch. 2, Sects. 5, 6]. In reality we repeat all proofs adding some changes which 
imply that the H\"older coefficient can be calculated explicitly (see 
Remark 5.5). Comparing to the original proofs from [LSU] we add some 
explanation which make them easier for reading.

\section{2. Notation and auxiliary results}

To prove the H\"older continuity of solutions to problem (1.6) we need

\proclaim Lemma 2.1. (see [Z1]). 
Let $u_0\in L_\infty(\Omega)$. Then the estimate holds
$$
\|u\|_{L_\infty(\Omega^T)}\le\|u_0\|_{L_\infty(\Omega)}\equiv M,\quad 
T<\infty.
\leqno(2.1)
$$

\Proof 
Let $A_k(t)=\{x\in\Omega:\ u(x,t)>k\}$, $u^{(k)}(x,t)=\max\{u-k,0\}$. 
Multiplying $(1.6)_1$ by $u^{(k)}$ and integrating over $\Omega$ yields
$$\eqal{
&{1\over2}{d\over dt}\intop_\Omega|u^{(k)}|^2dx+\nu\intop_\Omega
|\nabla u^{(k)}|^2dx-\nu\intop_{-a}^a|u^{(k)}|_{r=R}|^2dz\cr
&\quad+2\nu\intop_\Omega u_{,r}u^{(k)}drdz=0.\cr}
$$
The last term in the above equality equals
$$
\nu\intop_\Omega(|u^{(k)}|^2)_{,r}drdz=\nu\intop_{-a}^a|u^{(k)}|_{r=R}|^2dz,
$$
because $u^{(k)}|_{r=0}=0$. Otherwise the condition (1.7) implies 
contradiction for $u>k$. Hence, we derive
$$
{1\over2}{d\over dt}\intop_\Omega|u^{(k)}|^2dx+\nu\intop_\Omega
|\nabla u^{(k)}|^2dx=0.
$$
Assuming that $u_0<k$ we obtain (2.1). This concludes the proof.

\proclaim Lemma 2.2. (see [Z2]). 
Assume that $v_0\in L_2(\Omega)$. Then there exists a weak solution to (1.1) 
such that $v\in V_2^0(\Omega^T)=\{w:\ \|w\|_{L_\infty(0,T;L_2(\Omega))}+
\|\nabla w\|_{L_2(\Omega^T)}<\infty\}$ and
$$
\|v\|_{V_2^0(\Omega^T)}\le c\|v_0\|_{L_2(\Omega)}.
\leqno(2.2)
$$

\noindent
A proof of the H\"older regularity of solutions to (1.6) we begin with 
recalling some auxiliary lemmas from [LSU, Ch. 2].

\noindent
Now we introduce some notation. Let us introduce the cylinder
$$\eqal{
&Q(\varrho_0,\tau_0)=\{(x,t):\ x\in B_{\varrho_0}(x_0),\ t_0-\tau_0<t<t_0\},\cr
&B_{\varrho_0}(x_0)=\{x:\ |x-x_0|<\varrho_0\}.\cr}
$$
Moreover, we need the following cylinder with the same top
$$
Q(\varrho-\sigma_1\varrho,\tau-\sigma_2\tau)=\{(x,t):\ 
x\in B_{\varrho-\sigma_1\varrho}(x_0),\ t_0-(1-\sigma_2)\tau<t<t_0\},
$$
where $\sigma_1,\sigma_2\in[0,1)$, 
${\varrho_0\over2}\le\varrho-\sigma_1\varrho\le\varrho\le\varrho_0$, 
${\tau_0\over2}\le\tau-\sigma_2\tau\le\tau\le\tau_0$.

\noindent
Let the function $u=u(x,t)$ be given. Then we introduce
$$\eqal{
&A_{k,\varrho}(t)=\{x\in B_\varrho(x_0):\ u(x,t)>k\},\cr
&\mu(k,\varrho,\tau)=\intop_{t_0-\tau}^{t_0}{\rm meas}^{r\over q}
A_{k,\varrho}(t)dt,\cr}
$$
where ${2\over r}+{3\over q}={3\over2}$. Changing variables $x-x_0=\varrho_0\tilde x$, $t-t_0=\varrho_0^2\tilde t$ 
and denoting $\tilde\varrho=\varrho\varrho_0^{-1}$, 
$\tilde\tau=\tau_0\varrho_0^{-2}$, we introduce the cylinders
$$
Q'(\tilde\varrho,\tilde\tau)=\{(\tilde x,\tilde t):\ |\tilde x|<\tilde\varrho,
\ -\tilde\tau<\tilde t<0\}.
$$
In view of the above transformations we have
$$
{1\over2}\le\tilde\varrho-\sigma_1\tilde\varrho\le\tilde\varrho\le1,\quad
{\theta\over2}\le\tilde\tau-\sigma_2\tilde\tau\le\tilde\tau\le\theta
\equiv\tau_0\varrho_0^{-2}.
$$
Then cylinder $Q(\varrho_0,\tau_0)$ is transformed to the cylinder
$$
Q'(1,\theta)=\{(\tilde x,\tilde t):\ |\tilde x|<1,\ -\theta<\tilde t<0\}.
$$
Assume that the inequality is satisfied
$$\eqal{
&\|u^{(k)}\|_{V_2^0(Q'(\varrho-\sigma_1\varrho,\tau-\sigma_2\tau))}^2\cr
&\le\gamma\{[(\sigma_1\varrho)^{-1}+(\sigma_2\tau)^{-2}]
\|u^{(k)}\|_{L_2(Q'(\varrho,\tau))}^2+k^2\mu^{2(1+\varkappa)\over r}
(k,\varrho,\tau)\},\cr}
\leqno(2.3)
$$
where $u^{(k)}=\max\{u-k,0\}$ and 
$Q'(\varrho,\tau)=\{|x|<\varrho,\ -\tau<t<0\}$.

\noindent
Moreover, we assume
$$
{1\over2}\le\varrho-\sigma_1\varrho<\varrho\le1,\quad
{1\over2}\theta\le\tau-\sigma_2\tau<\tau\le\theta.
$$
Let us take the sequence of cylinders
$$
Q_h=Q'(\varrho_h,\tau_h)\quad {\rm with}\quad
\varrho_h={1\over2}+{1\over2^{h+2}},\quad 
\tau_h={1\over2}\theta+{\theta\over2^{h+2}},
$$
$h=0,1,\dots,\infty$ and the sequence of increasing levels 
$k_h=M+N\left(1-{1\over2^h}\right)$, $h=0,1,\dots$, where $M$ and $N$ are 
some positive numbers.

\proclaim Lemma 2.3. (see Lemma 6.1 from [LSU, Ch. 2, Sect. 6]) 
Let the function $u(x,t)$ satisfy (2.3) for all $k\in[M,M+N]$. Then the 
quantities
$$
y_h=N^{-2}\intop_{-\tau_h}^0\intop_{A_{k_h,\varrho_h}(t)}(u-k_h)^2dxdt
$$
and
$$
z_h=\mu^{2\over r}(k_h,\varrho_h,\tau_h)\equiv\mu_h^{2/r},\quad h=0,1,\dots,
\quad r\in\bigg[{4\over3},\infty\bigg],
$$
with defined above $k_h$, $\varrho_h$ and $\tau_h$ satisfy the system of 
recurrent inequalities
$$
y_{h+1}\le\gamma_12^{2h}\bigg[2^{2h}y_h^{1+\delta}+z_h^{1+\varkappa}y_h^\delta
\bigg({M\over N}+1\bigg)^2\bigg],
\leqno(2.4)
$$
$$
z_{h+1}\le\gamma_12^{2h}\bigg[2^{2h}y_h+z_h^{1+\varkappa}
\bigg({M\over N}+1\bigg)^2\bigg],
\leqno(2.5)
$$
where $\delta={2\over n+2}$ and 
$\gamma_1=4\beta^2(\gamma+1)\big(2^8+{2^{3-h}\over\theta}\big)$.

\Proof 
Let $\zeta_h(|x|)$ be a sequence of cut off continuous functions such that
$$
\zeta_h(|x|)=\cases{
1& for $|x|\le\varrho_{h+1}$\cr
0& for $|x|\ge\bar\varrho_h={1\over2}(\varrho_h+\varrho_{h+1})$\cr}
$$
and linear on the interval $|x|\in[\varrho_{h+1},\bar\varrho_h]$, so that
$|\zeta_{hx}(|x|)|\le2^{h+4}$.

\noindent
Let
$$
\lambda_h=\intop_{-\tau_{h+1}}^0{\rm meas}A_{k_{h+1},\bar\varrho_h}(t)dt,
$$
where $A_{k,\varrho}(t)=\{x\in B_\varrho:\ u(x,t)>k\}$. Then
$$\eqal{
&y_{h+1}=N^{-2}\intop_{-\tau_{h+1}}^0\intop_{A_{k_{h+1},\varrho_{h+1}}(t)}
(u-k_{h+1})^2dxdt\cr
&\le N^{-2}\intop_{-\tau_{h+1}}^0\intop_{A_{k_{h+1},\bar\varrho_h}(t)}
(u-k_{h+1})^2\zeta_h^2dxdt\equiv I_0.\cr}
\leqno(2.6)
$$
Employing $\bar\varrho_h>\varrho_{h+1}$ we have
$$
I_0\le\beta^2N^{-2}\lambda_h^{2\over n+2}
\|u^{(k_{h+1})}\zeta_h\|_{V_2^0(Q'(\bar\varrho_h,\tau_{h+1}))}\equiv I,
$$
where we used the inequality
$$
\|u\|_{L_2(\Omega^T)}\le\beta({\rm meas}Q_0)^{1\over n+2}
\|u\|_{V_2^0(\Omega^T)},
\leqno(2.7)
$$
where $Q_0=\{(x,t):\ |u(x,t)|>0\}$ and $\beta$ is a constant from the imbedding
$\|u\|_{L_r(0,T;L_q(\Omega))}\le\beta\|u\|_{V_2^0(\Omega^T)}$, where 
${2\over r}+{3\over q}={3\over2}$. Continuing,
$$\eqal{
I&\le\beta^2N^{-2}\lambda_h^{2\over n+2}\bigg\{\esssup_{-\tau_{h+1}\le t\le0}
\intop_{A_{k_{h+1},\bar\varrho_h}(t)}(u-k_{h+1})^2\zeta_h^2dx\cr
&\quad+2\intop_{-\tau_{h+1}}^0\intop_{A_{k_{h+1},\bar\varrho}(t)}
[u_x^2\zeta_h^2+(u-k_{h+1})^2\zeta_{hx}^2]dxdt\bigg\}\cr
&\le2\beta^2\lambda_h^\delta[N^{-2}
\|u^{(k_{h+1})}\|_{V_2^0(Q'(\bar\varrho_h,\tau_{h+1}))}^2+4^{h+4}y_h],\cr}
$$
because $k_h<k_{h+1}$ implies that $y_h\ge y_{h+1}$. Now we calculate
$$\eqal{
y_h&=N^{-2}\intop_{-\tau_h}^0\intop_{A_{k_h,\varrho_h}}(u-k_h)^2dxdt\ge\cr
&\quad N^{-2}\intop_{-\tau_h}^0dt\intop_{A_{k_{h+1},
\varrho_h}(t)}(u-k_h)^2dx\ge J.\cr}
\leqno(2.8)
$$
Because $k_{h+1}>k_h$, we have
$$\eqal{
J&\ge(k_{h+1}-k_h)^2N^{-2}\intop_{-\tau_h}^0dt\intop_{A_{k_{h+1},\varrho_h}}dx
\ge\cr
&\quad(k_{h+1}-k_h)^2N^{-2}\intop_{-\tau_h}^0dt\intop_{A_{k_{h+1},
\bar\varrho_h}}dx=(k_{h+1}-k_h)^2N^{-2}\lambda_h,\cr}
$$
because $\varrho_{h+1}<\varrho_h$ so 
$\varrho_{h+1}<\bar\varrho_h={1\over2}(\varrho_h+\varrho_{h+1})<\varrho_h$.

\noindent
Now we estimate the r.h.s. of (2.6) by using (2.3) and (2.8). We apply 
inequality (2.3) for $u^{(k+1)}$ and for two cylinders: 
$Q'(\bar\varrho_h,\tau_{h+1})$ and $Q'(\varrho_h,\tau_h)$.
First we calculate $(\sigma_1\varrho)^{-2}$ and $(\sigma_2\tau)^{-1}$.
We have $\sigma_1\varrho=\varrho_h-\bar\varrho_h={1\over2^{h+4}}$ 
so $(\sigma_1\varrho)^{-2}=4^{h+4}$ and 
$(\sigma_2\tau)^{-1}={2^{h+3}\over\theta}$. Hence we have
$$\eqal{
y_{h+1}&\leq2\beta^2[(k_{h+1}-k_h)^{-1}N]^{2\delta}y_h^\delta\cr
&\quad\cdot\bigg\{\gamma\bigg[\bigg(4^{h+4}+{1\over\theta}
2^{h+3}\bigg)N^{-2}\|u^{(k_{h+1})}\|_{L_2(Q'(\varrho_h,\tau_h))}^2\cr
&\quad+k_{h+1}^2N^{-2}\mu_h^{{2\over r}(1+\varkappa)}\bigg]+
4^{h+4}y_h\bigg\}\cr
&\le2\beta^2[(k_{h+1}-k_h)^{-1}N]^{2\delta}[2^{2h}a_1y_h^{1+\delta}+\gamma
k_{h+1}^2N^{-2}z_h^{1+\varkappa}y_h^\delta]\equiv I,\cr}
\leqno(2.9)
$$
where $a_1=(\gamma+1)\big(2^8+{2^{3-h}\over\theta}\big)$. Using that
$$
(k_{h+1}-k_h)^{-1}={2^{h+1}\over N}\quad {\rm so}\quad
(k_{h+1}-k_h)^{-1}N=2^{h+1},
$$
and
$$
k_{h+1}^2N^{-2}=\bigg({M+N\left(1-{1\over2^{h+1}}\right)\over N}\bigg)^2\le
\bigg({M\over N}+1\bigg)^2.
$$
Then (2.9) implies
$$
I\le4\beta^22^h\bigg[2^{2h}a_1y_n^{1+\delta}+\gamma\bigg({M\over N}+1\bigg)^2
z_h^{1+\varkappa}y_h^\delta\bigg],
$$
where we used that $2\delta\le1$. From the estimate inequality (2.4) follows.

\noindent
Next we show (2.5). For this aim we consider
$$\eqal{
&(k_{h+1}-k_h)^2z_{h+1}=(k_{h+1}-k_h)^2\mu_{h+1}^{2\over r}\cr
&\le\|u^{(k_h)}\zeta_h\|_{L_{q,r}(Q'(\bar\varrho_h,\tau_{h+1}))}^2\le\beta^2
\|u^{(k_h)}\zeta_h\|_{V_2^0(Q'(\bar\varrho_h,\tau_{h+1}))}^2.\cr}
$$
Next, using (2.3) we calculate
$$\eqal{
&(k_{h+1}-k_h)^2z_{h+1}\le2\beta^24^{h+4}N^2y_h\cr
&\quad+2\beta^2\gamma\big\{[(\bar\varrho_h-\varrho_h)^{-2}+
(\tau_{h+1}-\tau_h)^{-1}]\|u^{(k_h)}\|_{L_2(Q'(\varrho_h,\tau_h))}^2\cr
&\quad+k_h^2\mu_h^{2(1+\varkappa)\over r}\big\}\cr
&\le2\beta^2\gamma\big\{[4^{h+4}+\theta^{-1}2^{h+3}]N^2y_h+k_h^2
z_h^{1+\varkappa}\big\}.\cr}
\leqno(2.10)
$$
Since $(k_{h+1}-k_h)^{-2}={2^{2(h+1)}\over N^2}$ we obtain from (2.10) 
inequality (2.5).

\noindent
This concludes the proof.

\proclaim Lemma 2.4. (see Lemma 5.7 from [LSU, Ch. 2, Sect. 5]) 
Assume that the numbers $y_h$, $z_h$, $h=0,1,\dots$, satisfy the system of 
inequalities
$$\eqal{
&y_{h+1}\le c_0b^h(y_h^{1+\delta}+z_h^{1+\varepsilon}y_h^\delta),\cr
&z_{h+1}\le c_0b^h(y_h+z_h^{1+\varepsilon}),\cr}
\leqno(2.11)
$$
where $c_0$, $b$, $\varepsilon$, $\delta$ are fixed positive constants with 
$b\ge1$.\\
Then
$$
y_h\le\lambda b^{-{h\over d}},\quad 
z_h\le(\lambda b^{-{h\over d}})^{1\over1+\varepsilon}
\leqno(2.12)
$$
where
$$\eqal{
&d=\min\left\{\delta,{\varepsilon\over1+\varepsilon}\right\}\cr
&\lambda=\min\left\{(2c_0)^{-1/\delta}b^{-1/\delta d};
(2c_0)^{-{1+\varepsilon\over\varepsilon}}b^{-{1\over\varepsilon d}}\right\}\cr}
$$
and if
$$
y_0\le\lambda,\quad z_0\le\lambda^{1\over1+\varepsilon}.
\leqno(2.13)
$$

\Proof 
For $h=0$ (2.12) holds. Assume that (2.12) is satisfied for $y_h$ and $z_h$. 
Then in view of (2.11) we obtain
$$\eqal{
&y_{h+1}\le2c_0b^h(\lambda b^{-{h\over d}})^{1+\delta}=2c_0\lambda^{1+\delta}
b^{h\left(1-{1+\delta\over d}\right)}\equiv I_1,\cr
&z_{h+1}\le2c_0\lambda b^{h\left(1-{1\over d}\right)}\equiv I_2.\cr}
$$
Using that $\lambda\le(2c_0)^{-1/\delta}b^{-1/\delta d}$ we obtain
$$
I_1\le\lambda b^{-{h+1\over d}}b^{h(1-\delta/d)}\le\lambda b^{-{h+1\over d}}
$$
because $d\le\delta$.

\noindent
Using that 
$\lambda\le(2c_0)^{-{1+\varepsilon\over\varepsilon}}
b^{-{1\over\varepsilon d}}$ we derive
$$
I_2\le(\lambda b^{-{h\over d}})^{1\over1+\varepsilon}
b^{h-{\varepsilon\over1+\varepsilon}{h\over d}}\le
(\lambda b^{-{h\over d}})^{1\over1+\varepsilon}
$$
because $d\le{\varepsilon\over1+\varepsilon}$. This concludes the proof.

\proclaim Lemma 2.5. (see Lemma 5.8 from [LSU, Ch. 2, Sect. 5]) 
Let $u=u(x,t)$ be a measurable and bounded in cylinder 
$Q_{\varrho_0}=Q(\varrho_0,\theta_0\varrho_0^2)$. Let\break 
$Q_\varrho=Q(\varrho,\theta_0\varrho^2)$ and $Q_{b\varrho}$, $b>1$ 
be cylinders with the same top and the axis as $Q_{\varrho_0}$.
Moreover, $\varrho\le b^{-1}\varrho_0$ and $u=u(x,t)$ satisfies the relations: 
either
$$
osc\{u;Q_\varrho\}\le c_1\varrho^\delta
\leqno(2.14)
$$
or
$$
osc\{u;Q_\varrho\}\le\eta osc\{u;Q_{b\varrho}\}
\leqno(2.15)
$$
with some constants $c_1$, $\delta\le1$ and $\eta<1$.\\
Then for $\varrho\le\varrho_0$ it holds
$$
osc\{u;Q_\varrho\}\le c\varrho_0^{-\alpha}\varrho^\alpha,
\leqno(2.16)
$$
where $\alpha=\min\{-\ln_b\eta;\delta\}$, 
$c=b^\alpha\max\{\omega;c_1\varrho_0^\delta\}$, 
$\omega_0=osc\{u;Q_{\varrho_0}\}$.

\Proof 
Take the sequence of cylinders $Q_{\varrho_k}$, $\varrho_k=b^{-k}\varrho_0$, 
$k=0,1,2,\dots$, described above. Let $\omega_k=osc\{u,Q_{\varrho_k}\}$. 
From (2.14) and (2.15) we have
$$
\omega_k\le\max\{c_1\varrho_k^\delta;\eta\omega_{k-1}\},\quad k=1,2,\dots,
\leqno(2.17)
$$
where
$$
\omega_0\le cb^{-\alpha}.
$$
Hence for $y_k=b^{k\alpha}\omega_k$, $k=1,2,\dots,$ we obtain the estimate
$$
y_k\le\max\{b^{k\alpha}c_1\varrho_k^\delta;b^{k\alpha}\eta\omega_{k-1}\}=
\max\{c_1b^{k(\alpha-\delta)}\varrho_0^\delta;b^\alpha\eta y_{k-1}\}\equiv I.
$$
From assumptions of the lemma we have
$$
\eta\le b^{-\alpha},\quad \alpha<\delta,\quad c\ge b^\alpha c_1
\varrho_0^\delta.
$$
Hence
$$
I\le\max\{c_1\varrho_0^\delta;y_{k-1}\}\le\max\{cb^{-\alpha};y_{k-1}\}
$$
and
$$
y_0=\omega_0\le cb^{-\alpha}.
$$
From the above inequalities, we deduce that
$$
y_k\le cb^{-\alpha}.
$$
Hence $\omega_k\le b^{-k\alpha}y_k\le cb^{-k\alpha}b^{-\alpha}=cb^{-\alpha}
\left({\varrho_k\over\varrho_0}\right)^\alpha$.

\noindent
Let $Q_\varrho$ be a cylinder with arbitrary $\varrho\le\varrho_0$.
Then there exists $k\ge1$ such that $\varrho_k\le\varrho\le\varrho_{k-1}$. 
Therefore
$$
osc\{u;Q_\varrho\}\le osc\{u;Q_{\varrho_{k-1}}\}\le cb^{-\alpha}
\varrho_0^{-\alpha}\varrho_{k-1}^\alpha\le c\varrho_0^{-\alpha}\varrho^\alpha.
$$
Hence the lemma is proved.

\section{3. Existence in $B_2(Q_T,M,\gamma,r,\delta,\varkappa)$}

\Def{3.1.} 
We say that $u\in V_2^0(Q_T)$ such that $\esssup_{Q_T}|u|\le M$ belongs to 
$B_2(Q_T,M,\gamma,r,\delta,\varkappa)$ if the function 
$\omega(x,t)=\pm u(x,t)$ satisfies the inequalities
$$\eqal{
&\max_{t_0\le t\le t_0+\tau}
\|\omega^{(k)}(x,t)\|_{L_2(B_{\varrho-\sigma_1\varrho})}^2\le
\|\omega^{(k)}(x,t_0)\|_{L_2(B_{\varrho-\sigma_1\varrho})}^2\cr
&\quad+\gamma[(\sigma_1\varrho)^{-2}\|\omega^{(k)}\|_{L_2(Q(\varrho,\tau))}^2+
\mu^{{2\over r}(1+\varkappa)}(k,\varrho,\tau)]\cr}
\leqno(3.1)
$$
and
$$\eqal{
&\|\omega^{(k)}\|_{V_2^0(Q(\varrho-\sigma_1\varrho,\tau-\sigma_2\tau))}^2\cr
&\le\gamma[((\sigma_1\varrho)^{-2}+(\sigma_2\tau)^{-1})
\|\omega^{(k)}\|_{L_2(Q(\varrho,\tau))}^2+
\mu^{{2\over r}(1+\varkappa)}(k,\varrho,\tau)],\cr}
\leqno(3.2)
$$
where $\omega^{(k)}(x,t)=\max\{\omega(x,t)-k;0\}$, 
$Q(\varrho,\tau)=B_\varrho(x_0)\times(t_0,t_0+\tau)=\{(x,t)\in Q_T:\ 
|x-x_0|<\varrho,\ t_0<t<t_0+\tau\}$, $Q_T=\Omega\times(0,T)$, 
$\varrho,\tau$ are arbitrary positive numbers, 
$\sigma_1,\sigma_2$ -- arbitrary numbers from $(0,1)$,
$$\eqal{
&\mu(k,\varrho,\tau)=\intop_{t_0}^{t_0+\tau}{\rm meas}^{r\over q}
A_{k,\varrho}(t)dt,\cr
&A_{k,\varrho}(t)=\{x\in B_\varrho(x_0):\ \omega(x,t)>k\},\cr}
$$
the numbers $M$, $\gamma$, $\delta$, $\varkappa$ are arbitrary positive.
The numbers $r$, $q$ satisfy the relation
$$
{2\over r}+{3\over q}={3\over2}.
\leqno(3.3)
$$
Finally, $k$ is a positive number satisfying the condition
$$
\esssup_{Q(\varrho,\tau)}\omega(x,t)-k\le\delta.
\leqno(3.4)
$$

\proclaim Lemma 3.2. 
Assume that $u$ satisfies (1.6) and $|u|\le M$. Assume that 
$v\in L_{r'}(0,T;L_{q'}(\Omega))$ with 
${3\over q'}+{2\over r'}=1-{3\over2}\varkappa$, $\varkappa>0$ and 
$\divv v=0$. Then 
$u\in{\cal B}_2(Q_T\cap Q(\varrho,\tau),M,\gamma,r,\delta,\varkappa)$ where 
$k,\delta$ satisfy (3.4). The important case is $\varkappa=1/3$, for 
$r'=q'=10$.

\Proof 
For solutions to problem (1.6) we have the estimate
$$
\|u\|_{L_\infty(\Omega^T)}\le c\|u(0)\|_{L_\infty(\Omega)}.
\leqno(3.5)
$$
Let $u^{(k)}(x,t)=\max\{u(x,t)-k,0\}$, 
$Q(\varrho,\tau)=B_\varrho\times(t_0,t_0+\tau)=\{(x,t):\ |x-x_0|<\varrho,\ 
t_0<t<t_0+\tau\}$ is arbitrary cylinder located in $\R^3\times(0,T)$, 
$\varrho,\tau$ -- arbitrary numbers from $(0,1)$.
$$\eqal{
&A_{k,\varrho}(t)=\{x\in B_\varrho:\ u(x,t)>k\}\cr
&\mu(k,\varrho,\tau)=\intop_{t_0}^{t_0+\tau}{\rm meas}^{r\over q}
A_{k,\varrho}(t)dt\cr}
$$
where ${2\over r}+{3\over q}={3\over2}$. 
Since $k$ is an arbitrary number we can assume that
$$
\esssup_{Q(\varrho,\tau)}u(x,t)-k\le\delta.
\leqno(3.6)
$$
Multiplying $(1.6)_1$ by $u^{(k)}\zeta^2$, where $\zeta=\zeta(x,t)$ is 
a smooth function with support in $Q(\varrho,\tau)$, and integrating the 
result over $Q(\varrho,\tau)$ we obtain
$$
\intop_{Q(\varrho,\tau)}\bigg(u_t-v\cdot\nabla u-\nu\Delta u+\nu
{u_{,r}\over r}\bigg)u^{(k)}\zeta^2dxdt=0
\leqno(3.7)
$$
Continuing calculations in (3.7) and using that $\zeta=\zeta(x,t)$ as 
a function of the space variable $x$ has a compact support in $B_R(x_0)$ yields
$$\eqal{
&\intop_{Q(\varrho,\tau)}\bigg[{1\over2}\partial_t|u^{(k)}|^2\zeta^2+
{1\over2}v\cdot\nabla|u^{(k)}|^2\zeta^2\cr
&\quad+\nu\nabla(u^{(k)})\nabla(u^{(k)}\zeta^2)+{\nu\over2}
{(|u^{(k)}|^2)_{,r}\over r}\zeta^2\bigg]dxdt=0.\cr}
$$
Continuing we have
$$\eqal{
&\intop_{Q(\varrho,\tau)}\bigg[{1\over2}\partial_t(|u^{(k)}|^2\zeta^2)-
|u^{(k)}|^2\zeta\zeta_t+{1\over2}v\cdot\nabla(|u^{(k)}|^2\zeta^2)\cr
&\quad-v\cdot\nabla\zeta\zeta|u^{(k)}|^2+\nu\nabla u^{(k)}\nabla
(u^{(k)}\zeta)\zeta+\nu\nabla u^{(k)}u^{(k)}\zeta\nabla\zeta\bigg]dxdt\cr
&\quad+{\nu\over2}\intop_{Q(\varrho,\tau)}|u^{(k)}|_{,r}^2\zeta^2drdzdt=0.\cr}
$$
Continuing,
$$\eqal{
&{1\over2}\intop_{B_\varrho}|u^{(k)}(x,t_0+\tau)|^2\zeta^2(x,t_0+\tau)dx+
\nu\intop_{Q(\varrho,\tau)}|\nabla(u^{(k)}\zeta)|^2dxdt\cr
&\le{1\over2}\intop_{B_\varrho}|u^{(k)}(x,t_0)|^2\zeta^2(x,t_0)dx\cr
&\quad+\intop_{Q(\varrho,\tau)}|u^{(k)}|^2\zeta|\zeta_t|dxdt+
\intop_{Q(\varrho,\tau)}|v\cdot\nabla\zeta\zeta|\,|u^{(k)}|^2dxdt\cr
&\quad+\nu\bigg|\intop_{Q(\varrho,\tau)}u^{(k)}\nabla\zeta\nabla
(u^{(k)}\zeta)dxdt\bigg|\cr
&\quad+\nu\bigg|\intop_{Q(\varrho,\tau)}\nabla u^{(k)}u^{(k)}\nabla\zeta
\zeta dxdt\bigg|+\nu\intop_{Q(\varrho,\tau)}{|u^{(k)}|^2\over r}\zeta
|\zeta_{,r}|dxdt\cr
&\quad+{\nu\over2}\intop_{Q'(\varrho,\tau)}|u^{(k)}|^2\zeta^2|_{r=0}dzdt,\cr}
\leqno(3.8)
$$
where $Q'(\varrho,\tau)=\{(x,t)\in Q(\varrho,\tau):\ r=0\}$.

\noindent
The third term on the r.h.s. of (3.8) is estimated by
$$
\intop_{Q(\varrho,\tau)}|v|\,|\nabla\zeta|\zeta|u^{(k)}|^2dxdt,
$$
the fourth by
$$
{\varepsilon_1\over2}\intop_{Q(\varrho,\tau)}|\nabla(u^{(k)}\zeta)|^2dxdt+
{\nu^2\over2\varepsilon_1}\intop_{Q(\varrho,\tau)}|u^{(k)}|^2
|\nabla\zeta|^2dxdt,
$$
the fifth is equal
$$
{\nu\over2}\intop_{Q(\varrho,\tau)}\nabla|u^{(k)}|^2\nabla\zeta\zeta dxdt=
-{\nu\over2}\intop_{Q(\varrho,\tau)}|u^{(k)}|^2(\Delta\zeta\zeta+
|\nabla\zeta|^2)dxdt.
$$
Next, the last but one is equal
$$
\nu\intop_{Q(\varrho,\tau)}|u^{(k)}|^2{1\over r}\zeta|\zeta_{,r}|dxdt
$$
The last term on the  r.h.s. of (3.8) disappears because $u^{(k)}|_{r=0}=0$. 
Otherwise we have a contradiction with (1.7) for $u>k$.

\noindent
Setting $\varepsilon_1=\nu$ we obtain from (3.8) the inequality
$$\eqal{
&{1\over2}\intop_{B_\varrho}|u^{(k)}(x,t_0+\tau)|^2\zeta^2(x,t_0+\tau)dx+
{\nu\over2}\intop_{Q(\varrho,\tau)}|\nabla(u^{(k)}\zeta)|^2dxdt\cr
&\le{1\over2}\intop_{B_\varrho}|u^{(k)}(x,t_0)|^2\zeta^2(x,t_0)dx\cr
&\quad+\intop_{Q(\varrho,\tau)}|u^{(k)}|^2\bigg(\zeta|\zeta_t|+
\nu|\nabla\zeta|^2+{\nu\over2}|\Delta\zeta|\zeta+
{\nu\over r}\zeta|\zeta_{,r}|\bigg)dxdt\cr
&\quad+\intop_{t_0}^{t_0+\tau}\intop_{A_{k,\varrho}(t)}|v|\,
|\nabla\zeta|\zeta|u^{(k)}|^2 dxdt.\cr}
\leqno(3.9)
$$
Now estimate the last term by
$$\eqal{
I&\equiv\intop_{t_0}^{t_0+\tau}\intop_{A_{k,\varrho}(t)}|v|^2|u^{(k)}|^2dxdt
+\intop_{t_0}^{t_0+\tau}\intop_{A_{k,\varrho}(t)}\zeta^2
|\nabla\zeta|^2|u^{(k)}|^2dxdt,\cr}
$$
where the second integral can be estimated in the same way as the second 
integral on the r.h.s. of (3.9) and the first integral in $I$ by
$$\eqal{
&\delta^2\bigg[\intop_{t_0}^{t_0+\tau}\bigg(\intop_{A_{k,\varrho}(t)}
|v|^{2\lambda_1}dx\bigg)^{2\mu_1/2\lambda_1}\bigg]^{2/2\mu_1}\cr
&\quad\cdot\bigg(\intop_{t_0}^{t_0+\tau}
|\meas A_{k,\varrho}(t)|^{a\mu_2\over a\lambda_2}dt\bigg)^{2/a\mu_2}\equiv I',
\cr}
$$
where ${1\over\lambda_1}+{1\over\lambda_2}=1$, 
${1\over\mu_1}+{1\over\mu_2}=1$, $a>1$.

\noindent
Assuming that $\mu_2$ and $\lambda_2$ are such that 
${3\over a\lambda_2}+{2\over a\mu_2}={3\over2}$ so 
${3\over\lambda_2}+{2\over\mu_2}={3\over2}a$ we have that 
${3\over\lambda_1}+{2\over\mu_1}=5-{3\over2}a$. Setting $q=a\lambda_2$, 
$r=a\mu_2$ we need that $v\in L_{r'}(0,T;L_{q'}(\Omega))$, where 
$r'=2\mu_1$, $q'=2\lambda_1$ and
$$
{3\over q'}+{2\over r'}={5\over2}-{3\over4}a
$$
Setting $a=2(1+\varkappa)$, $\varkappa>0$ we obtain
$$
{3\over q'}+{2\over r'}=1-{3\over2}\varkappa.
$$
Then
$$
I'\le\delta^2\|v\|_{L_{r'}(0,T;L_{q'}(B_\varrho)}
\mu^{{2\over r}(1+\varkappa)}(k,\varrho,\tau)
$$
In the case $q'=r'=10$ we have that $\varkappa={1\over3}$.
Moreover, $\varkappa={1\over6}$ for ${3\over q'}+{2\over r'}={3\over4}$.

\noindent
Let us assume that $\zeta(x,t)=1$ for 
$(x,t)\in Q(\varrho-\sigma_1\varrho,\tau-\sigma_2\tau)$, where 
$\sigma_1,\sigma_2\in(0,1)$. Then (3.9) takes the form
$$\eqal{
&{c_2\over2}\|u^{(k)}\|_{V_2^0(Q(\varrho-\sigma_1\varrho,\tau-\sigma_2\tau))}^2
\cr
&\le{1\over2}\|u^{(k)}(\cdot,t_0)\|_{L_2(B_\varrho)}^2+\gamma
[(\sigma_2\tau)^{-1}+(\sigma_1\varrho)^{-2}]
\|u^{(k)}\|_{L_2(Q(\varrho,\tau))}^2\cr
&\quad+c_1\delta^2\mu
(k,\varrho,\tau)^{{2\over r}(1+\varkappa)},\cr}
\leqno(3.10)
$$
where $c_1=\|v\|_{L_{r'}(t_0,t_0+\tau;L_{q'}(B_\varrho))}$, 
$c_2=\min\{1,\nu\}$ and we used the inequality
$$
\zeta|\zeta_t|+\nu|\nabla\zeta|^2+{\nu\over2}|\Delta\zeta|^2+{\nu\over r}
\zeta|\zeta_{,r}|+\zeta^2|\nabla\zeta|^2\le
\gamma[(\sigma_2\tau)^{-1}+(\sigma_1\varrho)^{-2}].
\leqno(3.11)
$$
In view of (3.10) we see that (3.1) and (3.2) are satisfied, so by 
Definition 3.1 the lemma is proved.

\noindent
If we take $\zeta=\zeta(x)$. Then instead of (3.10) we obtain
$$\eqal{
&\max_{t_0\le t\le t_0+\tau}
\|u^{(k)}(x,t)\|_{L_2(B_{\varrho-\sigma_1\varrho})}^2\le
\|u^{(k)}(x,t_0)\|_{L_2(B_\varrho)}^2\cr
&\quad+\gamma(\sigma_1\varrho)^{-2}\|u^{(k)}\|_{L_2(Q(\varrho,\tau))}^2+
c_1\delta^2\mu(k,\varrho,\tau)^{{2\over r}(1+\varkappa)}.\cr}
\leqno(3.12)
$$

\proclaim Lemma 3.3.
Let function $u=u(x,t)$ satisfy (3.12) in cylinder 
$Q(\varrho,t)=B_\varrho(x_0)\times(t_0,t)$, $t\in[t_0,t_0+\theta\varrho^2]$ 
and meas $A_{k,\varrho}(t_0)\le{1\over p}$ meas 
$B_\varrho={1\over p}\varkappa_n\varrho^n$, $p>1$. Then for any 
$\xi\in\left({1\over\sqrt{p}},1\right)$ there exist positive numbers 
$\theta(\xi)$ and $b(\xi)$ such that if
$$
\delta\ge H=
\esssup_{{x\in B_\varrho}\atop{t_0\le t\le t_0+\theta(\xi)\varrho^2}}
u(x,t)-k>\varrho^{n\varkappa\over2}
\leqno(3.13)
$$
then
$$
{\rm meas}(B_\varrho\setminus A_{k+\xi H,\varrho}(t))\ge b(\xi)
\varkappa_n\varrho^n
\leqno(3.14)
$$
for $t\in[t_0,t_0+\theta(\xi)\varrho^2]$ and $b(\xi)$, $\theta(\xi)$ are 
defined by (3.20).

\Remark{3.4.} 
We consider the case $n=3$, $\varkappa_3={4\pi\over3}$, $\varkappa={1\over3}$.
\vskip6pt

\noindent
{\bf Proof of Lemma 3.3.}
In view of the assumptions of the lemma and (3.12) we have
$$\eqal{
&\intop_{A_{k,\varrho-\sigma_1\varrho}(t)}(u(x,t)-k)^2dx\le{1\over p}
H^2\varkappa_n\varrho^n\cr
&\quad+\gamma H^2\sigma_1^{-2}(t-t_0)\varkappa_n
\varrho^{n-2}+c_1\delta^2[(t-t_0)(\varkappa_n
\varrho^n)^{r\over q}]^{{2\over r}(1+\varkappa)}.\cr}
\leqno(3.15)
$$
From the other hand
$$
(\xi H)^2{\rm meas}A_{k+\xi H,\varrho-\sigma_1\varrho}(t)\le
\intop_{A_{k,\varrho-\sigma_1\varrho}(t)}(u-k)^2dx.
\leqno(3.16)
$$
For the reader convenience we show (3.16)
$$\eqal{
&\xi^2H^2{\rm meas}A_{k+\xi H,\varrho-\sigma_1\varrho}(t)=
\intop_{B_{\varrho-\sigma_1\varrho}\cap\{x:\ u(x,t)>k+\xi H\}}\xi^2H^2dx\cr
&\le\intop_{B_{\varrho-\sigma_1\varrho}\cap\{x:\ u(x,t)>k+\xi H\}}
|u(x,t)-k|^2dx\cr
&\le\intop_{B_{\varrho-\sigma_1\varrho}\cap\{x:\ u(x,t)>k\}}
|u(x,t)-k|^2dx=\intop_{A_{k,\varrho-\sigma_1\varrho}}(u-k)^2dx,\cr}
$$
where we used the inequalities
$$\eqal{
&A_{k+k_0,\varrho}(t)=\{x\in B_\varrho:\ u(x,t)>k+k_0\},\quad k_0\ge0,\cr
&k_0^2{\rm meas}A_{k+k_0,\varrho}=\intop_{A_{k+k_0,\varrho}}k_0^2dx\le
\intop_{A_{k+k_0,\varrho}}(u-k)^2dx\cr
&\le\intop_{A_{k,\varrho}}(u-k)^2dx.\cr}
$$
Hence, in view of (3.15) and (3.16), we have
$$\eqal{
&{\rm meas}A_{k+\xi H,\varrho-\sigma_1\varrho}(t)\le{1\over p\xi^2}
\varkappa_n\varrho^n+[\gamma(\sigma_1\xi)^{-2}(t-t_0)\varrho^{-2}\cr
&\quad+c_1\delta^2(\xi H)^{-2}\varkappa_n^{{2\over q}(1+\varkappa)-1}
(t-t_0)^{{2\over r}(1+\varkappa)}\varrho^{n{2\over q}(1+\varkappa)-n}]
\varkappa_n\varrho^n.\cr}
\leqno(3.17)
$$
From (3.17) for $t\le t_0+\theta(\xi)\varrho^2$ and 
$H>\varrho^{n\varkappa\over2}$ we obtain
$$\eqal{
&{\rm meas}A_{k+\xi H,\varrho-\sigma_1\varrho}(t)\le\xi^{-2}\{p^{-1}+
\gamma\sigma_1^{-2}\theta(\xi)\cr
&\quad+c_1\delta^2\varkappa_n^{{2\over q}(1+\varkappa)-1}\theta
(\xi)^{{2\over r}(1+\varkappa)}
\varrho^{n{2\over q}(1+\varkappa)-n-n\varkappa+{4\over r}(1+\varkappa)}\}
\varkappa_n\varrho^n,\cr}
\leqno(3.18)
$$
where
$$
{2n\over q}(1+\varkappa)-n-n\varkappa+{4\over r}(1+\varkappa)=
{4\over r}+{2n\over q}-n+\bigg({4\over r}+{2n\over q}-n\bigg)\varkappa=0
$$
because $r$, $q$ satisfy the relation
$$
{4\over r}+{2n\over q}=n.
$$
Then (3.18) takes the form
$$\eqal{
&{\rm meas}A_{k+\xi H,\varrho-\sigma_1\varrho}(t)\le\xi^{-2}\{p^{-1}+\gamma
\sigma_1^{-2}\theta(\xi)\cr
&\quad+c_1\delta^2\varkappa_n^{{2\over q}(1+\varkappa)-1}\theta
(\xi)^{{2\over r}(1+\varkappa)}\}\varkappa_n\varrho^n.\cr}
\leqno(3.19)
$$
For any $\xi\in\left(\sqrt{1\over p},1\right)$ we choose $\sigma_1$, 
$\theta(\xi)$, $b(\xi)$ such that
$$
n\sigma_1+\xi^{-2}[p^{-1}+\gamma\sigma_1^{-2}\theta(\xi)+c_1\delta^2
\varkappa_n^{{2\over q}(1+\varkappa)-1}\theta(\xi)^{{2\over r}(1+\varkappa)}]
=1-b(\xi)<1.
\leqno(3.20)
$$
Then
$$\eqal{
{\rm meas}A_{k+\xi H,\varrho}(t)&\le{\rm meas}
A_{k+\xi H,\varrho-\sigma_1\varrho}(t)+{\rm meas}B_\varrho\setminus
B_{\varrho-\sigma_1\varrho}\cr
&\le(1-b(\xi)){\rm meas}B_\varrho,\cr}
\leqno(3.21)
$$
because
$$\eqal{
&A_{k,\varrho}(t)=\{x\in B_\varrho:\ u(x,t)>k\},\cr
&A_{k,\varrho-\sigma_1\varrho}(t)=\{x\in B_{\varrho-\sigma_1\varrho}:\ 
u(x,t)>k\},\cr}
$$
$$
{\rm meas}B_\varrho\setminus B_{\varrho-\sigma_1\varrho}=\varkappa_n\varrho^n-
\varkappa_n\varrho^n(1-\sigma_1)^n=\varkappa_n\varrho^n[1-(1-\sigma_1)^n]
\le\varkappa_n\varrho^nn\sigma_1.
$$
For $n=3$
$$
1-(1-\sigma_1)^3=1-(1-3\sigma_1+3\sigma_1^2-\sigma_1^3)=
3\sigma_1-\sigma_1^2(3-\sigma_1)<3\sigma_1.
$$
Hence, (3.21) implies that
$$
{\rm meas}(B_\varrho\setminus A_{k+\xi H,\varrho})\ge b(\xi){\rm meas}
B_\varrho.
$$
This proves Lemma 3.3
\kwadrat

Now we determine the quantities satisfying (3.20).
We choose $\sigma_1={1\over6n}$, $\xi={3\over4}$, $p=4$, $b(\xi)={1\over18}$.
Then $\xi\in\left({1\over2},1\right)$ and (3.20) takes the form
$$
{1\over6}+{16\over9}\cdot{1\over4}+{16\over9}[(6n)^2\gamma\theta+c_1\delta^2
\varkappa_n^{{2\over q}(1+\varkappa)-1}\theta^{{2\over r}(1+\varkappa)}]=
{17\over18}
\leqno(3.22)
$$
Hence
$$
(6n)^2\gamma\theta+c_1\delta^2\varkappa_n^{{2\over q}(1+\varkappa)-1}
\theta^{{2\over r}(1+\varkappa)}={9\over16}\cdot{1\over3}={3\over16}
\leqno(3.23)
$$
We have to recall that in our case $n=3$, $\varkappa=1/3$, 
${2\over r}+{3\over q}={3\over2}$, $\varkappa_3={4\pi\over3}$.

\noindent
Hence, in our case (3.23) takes the form
$$
9\cdot36\cdot\gamma\theta+c_1\delta^2\bigg({4\pi\over3}\bigg)^{-1/5}
\theta^{4/5}={3\over16}.
\leqno(3.24)
$$

\section{4. Proof of Lemma 7.2 from [LSU, Ch. 2, Sect. 7]}

Let 
$Q_\varrho=\{(x,t):\ |x-x_0|\le\varrho,\ t_0\le t\le t_0+\theta\varrho^2\}$. 
Let $Q_\varrho(k)=\{(x,t)\in Q_\varrho:\ u(x,t)>k\}$. Then we have 
the following inequalities
$$\eqal{
&(\varkappa_n\varrho^n)^{{2\over q}-{2\over r}}{\rm meas}^{2\over r}
Q_\varrho(k)\le\mu^{2\over r}(k,\varrho,\theta\varrho^2)\quad\cr
&\le(\theta\varrho^2)^{{2\over r}-{2\over q}}{\rm meas}^{2\over q}
Q_\varrho(k)\quad &{\rm for}\ \ r\le q,\cr
&(\varkappa_n\varrho^n)^{{2\over q}-{2\over r}}{\rm meas}^{2\over r}
Q_\varrho(k)\ge\mu^{2\over r}(k,\varrho,\theta\varrho^2)\cr
&\ge(\theta\varrho^2)^{{2\over r}-{2\over q}}{\rm meas}^{2\over q}
Q_\varrho(k)\quad &{\rm for}\ \ r\ge q.\cr}
\leqno(4.1)
$$
$$\eqal{
&{\rm meas}A_{k,\varrho}(t)=\{x\in B_\varrho(x_0):\ u(x,t)>k\},\cr
&\mu(k,\varrho,\theta\varrho^2)=\intop_{t_0}^{t_0+\theta\varrho^2}
{\rm meas}^{r\over q}A_{k,\varrho}(t)dt,\quad
{1\over r}+{n\over2q}={n\over4}.\cr}
$$
$$\eqal{
{\rm meas}Q_\varrho(k)&=\intop_{t_0}^{t_0+\theta\varrho^2}{\rm meas}
A_{k,\varrho}(t)dt\cr
&=\intop_{t_0}^{t_0+\theta\varrho^2}{\rm meas}^{r\over q}
A_{k,\varrho}(t){\rm meas}^{1-{r\over q}}A_{k,\varrho}(t)dt\cr}
$$
and we have that
$$\eqal{
&{\rm meas}^{1-{r\over q}}A_{k,\varrho}(t)\le{\rm meas}^{1-{r\over q}}
B_\varrho(x_0)\quad &{\rm for}\ \ r\le q,\cr
&{\rm meas}^{1-{r\over q}}A_{k,\varrho}(t)\ge{\rm meas}^{1-{r\over q}}
B_\varrho(x_0)\quad &{\rm for}\ \ r\ge q.\cr}
$$
Now we show the second inequality of $(4.1)_1$.
We calculate (the case $q\ge r$)
$$\eqal{
&\mu^{2\over r}(k,\varrho,\theta\varrho^2)=
\bigg[\intop_{t_0}^{t_0+\theta\varrho^2}A_{k,\varrho}^{r\over q}(t)
dt\bigg]^{2/r}\cr
&\le\bigg[\bigg(\intop_{t_0}^{t_0+\theta\varrho^2}1dt\bigg)^{1/\lambda_1}
\bigg(\intop_{t_0}^{t_0+\theta\varrho^2}
A_{k,\varrho}^{{r\over q}\lambda_2}(t)dt\bigg)^{1/\lambda_2}\bigg]^{2/r}
\equiv I_1,\cr}
$$
where ${1\over\lambda_1}+{1\over\lambda_2}=1$, ${r\over q}\lambda_2=1$, so 
$\lambda_2={q\over r}$, $\lambda_1={q\over q-r}$.

\noindent
Continuing, we have
$$\eqal{
I_1&=\bigg[(\theta\varrho^2)^{1-{r\over q}}\bigg(
\intop_{t_0}^{t_0+\theta\varrho^2}A_{k,\varrho}(t)
dt\bigg)^{r/q}\bigg]^{2\over r}
=(\theta\varrho^2)^{{2\over r}-{2\over q}}{\rm meas}^{2\over q}Q_\varrho(k).
\cr}
$$
Next we show the first inequality of $(4.1)_1$. We have 
$$\eqal{
&\meas^{2\over r}Q_\varrho(k)=\bigg(\intop_{t_0}^{t_0+\theta\varrho^2}
\meas^{r\over q}A_{k,\varrho}(t)\meas^{1-r/q}A_{k,\varrho}(t)dt\bigg)^{2/r}\cr
&\le(\varkappa_n\varrho^n)^{{2\over r}-{2\over q}}\mu^{2/r}
(k,\varrho,\theta_\varrho^2),\cr}
$$
where $r\le q$, so first inequality of $(4.1)_1$ follows.

\proclaim Lemma 4.1. (see Lemma 7.2 from [LSU, Ch. 2, Sect. 7]) 
Let $u=u(x,t)$ satisfy inequality (3.2). There exists $\theta_1>0$ such 
that for any cylinder $Q_{\varrho_0}\subset Q_T$ and any number
$$
k_0\ge\esssup_{Q_{\varrho_0}}u(x,t)-\delta
$$
the inequality
$$
{\rm meas}Q_{\varrho_0}(k_0)\le\theta_1\varrho_0^{n+2}
\leqno(4.2)
$$
implies
$$
{\rm meas}Q_{\varrho_0\over\sigma}\bigg(k_0+{H\over2}\bigg)=0,\quad \sigma>1,
\leqno(4.3)
$$
if
$$
H=\max_{Q_{\varrho_0}}u(x,t)-k_0\ge\varrho_0^{n\varkappa\over2}.
\leqno(4.4)
$$
The cylinders $Q_{\varrho_0/\sigma}$ and $Q_{\varrho_0}$ have the same top 
and the axis.

\Proof
Introducing new variables $\tilde x,\tilde t$ such that
$$\eqal{
&x-x_0=\varrho_0\tilde x,\quad t-t_0=\varrho_0^2\tilde t,\quad
|\tilde x|={1\over\varrho_0}|x-x_0|\le{\varrho\over\varrho_0}
\equiv\tilde\varrho,\cr
&|\tilde t|=\bigg|{1\over\varrho^2}(t-t_0)\bigg|\le{\tau\over\varrho_0^2}
\equiv\tilde\tau,\quad
\tilde\tau={\tau\over\varrho_0^2}={\theta^2\varrho^2\over\varrho_0^2}=
\theta^2\tilde\varrho,\quad
\tilde\varrho_0={\varrho_0\over\varrho_0}=1,\cr}
$$
the inequality (3.2) takes the form
$$\eqal{
&\|u^{(k)}\|_{V_2^0
(\tilde\varrho-\sigma_1\tilde\varrho,\tilde\tau-\sigma_2\tilde\tau)}^2\cr
&\le\gamma\{[(\sigma_1\tilde\varrho)^{-2}+(\sigma_2\tilde\tau)^{-1}]
\|u^{(k)}\|_{L_2(Q(\tilde\varrho,\tilde\tau))}^2
+\varrho_0^{n\varkappa}\mu^{{2\over r}(1+\varkappa)}
(k,\tilde\varrho,\tilde\tau)\},\cr}
\leqno(4.5)
$$
and $Q_{\tilde\varrho_0}=Q(1,\theta)=\{|\tilde x|<1,\ 0<\tilde\tau<\theta^2\}$.
Let $H$ be the number from (4.4). Then
$$
\varrho_0^{n\varkappa\over2}H^{-1}\le1.
\leqno(4.6)
$$
Introducing the new quantity $v={u\over H}$, using (4.6) we obtain from 
(4.5) the inequality
$$\eqal{
&\|v^{(\tilde k)}\|_{V_2^0(Q(\tilde\varrho-\sigma_1\tilde\varrho,
\tilde\tau-\sigma_2\tilde\tau))}\cr
&\le\gamma\big\{[(\sigma_1\tilde\varrho)^{-2}+(\sigma_2\tilde\tau)^{-1}]
\|v^{(\tilde k)}\|_{L_2(Q(\tilde\varrho,\tilde\tau))}^2+
\mu^{{2\over r}(1+\varkappa)}(\tilde k,\tilde\varrho,\tilde\tau)\big\},\cr}
\leqno(4.7)
$$
where $\tilde k={k\over H}$.
The function $\tilde v(\tilde x,\tilde t)$ satisfies the assumptions 
of Lemma~2.3. Moreover, instead of interval [M, M+N] we take 
$\left[\tilde k_0={k_0\over H},{k_0\over H}+{1\over2}\right]$, so we take 
$\tilde k\in\left[{k_0\over H},{k_0\over H}+{1\over2}\right]$. 
Then, in view of Lemma 2.3, the quantities
$$\eqal{
&y_h=4\intop_{-\tilde\tau_h}^0d\tilde t
\intop_{A_{\tilde k_h,\tilde\varrho_h}(\tilde t)}(v-\tilde k_h)^2d\tilde x,\cr
&z_h=\bigg(\intop_{-\tilde\tau_h}^0{\rm meas}^{r\over q}
A_{\tilde k_h,\tilde\varrho_h}(\tilde t)d\tilde t\bigg)^{2/r},\cr}
\leqno(4.8)
$$
satisfy inequalities (2.4), (2.5) with $M={k_0\over H}$, $N={1\over2}$.
Then in view of Lemma 2.4,
$$
y_h,z_h\to_{h\to\infty}0
$$
if
$$
y_0\le\min\{(2\gamma_1)^{-{n+2\over2}}(64)^{-{1\over\delta d}};\ 
(2\gamma_1)^{-{1+\varkappa\over\varkappa}}64^{-{1\over\varkappa d}}\}
\equiv\lambda,
\quad z_0\le\lambda^{1\over1+\varkappa},
$$
where $d=\min\left\{{2\over n+2},{\varkappa\over1+\varkappa}\right\}$, 
$\delta={2\over n+2}$ and 
$\gamma_1=4\beta^2(\gamma+1)\big(2^8+{2^3\over\theta}\big)$. In our case
$\delta=2/5,\quad \varkappa=1/3=\varepsilon$ so
$d=\min\{2/5,1/4\}=1/4$.

\noindent
Hence
$\lambda=\min\{(2\gamma_1)^{-{5\over2}}(64)^{-{35\over2}},(2\gamma_1)^{-4}
64^{-42}\}=(2\gamma_1)^{-4}64^{-42}$.

\noindent
In view of (4.2) we have
$$\eqal{
y_0&=4\intop_{-\theta}^0d\tilde t\intop_{A_{\tilde k_0,1}}(v-\tilde k_0)^2
d\tilde x\le\intop_{-\theta}^0{\rm meas}A_{\tilde k_0,1}(\tilde t)d\tilde t\cr
&=\varrho_0^{-n-2}{\rm meas}Q_{\varrho_0}(k_0)\le\theta_1,\cr
z_0&=\bigg(\intop_{-\theta}^0{\rm meas}^{r\over q}A_{k_0,1}(\tilde t)
\tilde dt\bigg)^{2/r}\le\cases{
\varkappa_n^{{2\over q}-{2\over r}}\theta_1^{2/r}& $r\ge q$,\cr
\theta^{{2\over r}-{2\over q}}\theta_1^{2/q}& $r\le q$.\cr}\cr}
$$
Hence
$$\eqal{
&\theta_1\le\min\big\{\lambda;\ \varkappa_n^{1-{r\over q}}
\lambda^{r\over2(1+\varkappa)}\big\}\quad &r\ge q,\cr
&\theta_1\le\min\big\{\lambda;\ \theta^{1-{q\over r}}
\lambda^{q\over2(1+\varkappa)}\big\}\quad &r\le q.\cr}
\leqno(4.9)
$$
For ${r\over2}+{3\over q}={3\over2}$, $\varkappa={1\over3}$ we have
${r\over2(1+\varkappa)}={3r\over8}$.
Hence for $\theta_1$ satisfying (4.9) we obtain that $y_h\to0$ for 
$h\to\infty$. This implies that
$$
\esssup_{Q\varrho_0/\sigma}u\le k_0+{H\over2}.
$$
This concludes the proof.

\Remark{4.2.} 
Let us calculate
$$
\gamma_1=2^6\beta^2[\gamma(2^5+\theta^{-1})+2^5],
$$
where $\beta$ is the constant from the imbedding
$$
\|u\|_{L_q(\Omega)}\le\beta\|u_x\|_{L_2(\Omega)}^\alpha
\|u\|_{L_2(\Omega)}^{1-\alpha}\quad \alpha={n\over2}-{n\over q}.
$$
Setting $\theta=2^{-5}$ we obtain $\gamma_1=2^{11}\beta^2(\gamma+1)$. Then \\
$\lambda=2^{-12\cdot7}[\beta^2(\gamma+1)]^{-7}64^{-42}=2^{-252}
[\beta^2(\gamma+1)]^{-7}$ and in view of (4.9)
$\theta_1\le2^{-252}[\beta^2(\gamma+1)]^{-7}$.

\section{5. Lemma 7.3 from [LSU, Ch. 2, Sect. 7]}

In this section the proof of Lemma 7.3 from [LSU, Ch. 2, Sect. 7] is changed 
in such direction that we can calculate the H\"older exponent $\alpha$ larger 
than the one in Theorem 7.1 in [LSU, Ch. 2, Sect. 7].

\noindent
Let 
$Q_{\sigma\varrho}=B_{\sigma\varrho}\times[t_0,t_0+2\sigma\theta\varrho^2]$,
$Q_\varrho=B_\varrho\times\left[t_0+{3\over2}\sigma\theta\varrho^2,
t_0+2\sigma\theta\varrho^2\right]$, 
$1<\sigma\le2$, have the same top and the axis.

Let $\mu_1=\max_{Q_{\sigma\varrho}}u$, $\mu_2=\min_{Q_{\sigma\varrho}}u$, 
$\omega=\mu_1-\mu_2$.

\proclaim Lemma 5.1. 
Let $u\in{\cal B}_2(\Omega^T,M,\gamma,r,\delta,\varkappa)$. Then for any 
$\theta_1>0$, $\sigma\in(1,2]$, there exists $s=s(\theta_1)>0$ such that either
$$
\omega=osc\{u,Q_{\sigma\varrho}\}\le2^s\varrho^{n\varkappa\over2}
\leqno(5.1)
$$
or
$$
{\rm meas}Q_\varrho\bigg(\mu_1-{\omega\over2^s}\bigg)\le\theta_1\varrho^{n+2}
\leqno(5.2)
$$
or
$$
{\rm meas}\bigg\{(x,t)\in Q_\varrho:\ u(x,t)<\mu_2+{\omega\over2^s}\bigg\}\le
\theta_1\varrho^{n+2}.
\leqno(5.3)
$$

\Proof 
Let $\omega\ge2^s\varrho^{n\varkappa\over2}$. Then we have either
$$
{\rm meas}A_{\mu_1-{\omega\over2},\varrho}
\bigg(t_0+{3\over2}\sigma\theta\varrho^2\bigg)\le{1\over2}\varkappa_n\varrho^n
\leqno(5.4)
$$
or
$$
{\rm meas}\bigg\{x\in B_\varrho:\ u
\bigg(x,t_0+{3\over2}\sigma\theta\varrho^2\bigg)<\mu_1-{\omega\over2}\bigg\}
\le{1\over2}\varkappa_n\varrho^n.
\leqno(5.5)
$$
The relation (5.4) means
$$
{\rm meas}\bigg\{x\in B_\varrho:\ 
u\bigg(x,t_0+{3\over2}\sigma\theta\varrho^2\bigg)>\mu_1-{\omega\over2}\bigg\}
\le{1\over2}\varkappa_n\varrho^n,
\leqno(5.6)
$$
where $\mu_1-{\omega\over2}={1\over2}(\mu_1+\mu_2)$. If (5.6) is valid, then
$$
{\rm meas}A_{\mu_1-{\omega\over2^r},\varrho}
\bigg(t_0+{3\over2}\sigma\theta\varrho^2\bigg)\le{1\over2}\varkappa_n\varrho^n,
\quad r\ge1.
\leqno(5.7)
$$
Hence
$$
{\rm meas}\bigg\{x\in B_\varrho:\ u
\bigg(x,t_0+{3\over2}\sigma\theta\varrho^2\bigg)>\mu_1-{\omega\over 2^r}\bigg\}
\le{1\over2}\varkappa_n\varrho^n
$$
because
$$
\meas A_{{\mu_1}-{\omega\over2^r},\varrho}
\bigg(t_0+{3\over2}\sigma\theta\varrho^2\bigg)\le\meas 
A_{\mu_1-{\omega\over2},\varrho}
\bigg(t_0+{3\over2}\sigma\theta\varrho^2\bigg)\le{1\over2}\varkappa_n\varrho^n,
\leqno(5.8)
$$
where $r=1,2,\dots,s-1$. Hence, condition (5.8) yields
$$
{\rm meas}\bigg\{x\in B_\varrho:\ u\bigg(x,t+{3\over2}
\sigma\theta\varrho^2\bigg)>\mu_1-{\omega\over2^{s-1}}\bigg\}\le
{1\over2}\varkappa_n\varrho^n.
\leqno(5.9)
$$
Now we proceed the calculations with $u=u(x,t)$. Let 
$\mu'_1=\max_{Q_\varrho}u$.

\noindent
If $\mu'_1\le\mu_1-{\omega\over2^s}$, then meas 
$Q_\varrho\left(\mu_1-{\omega\over2^s}\right)=0$ because
$$\eqal{
&Q_\varrho\bigg(\mu_1-{\omega\over2^s}\bigg)=\bigg\{(x,t)\in Q_\varrho:\ 
u(x,t)>\mu_1-{\omega\over2^s}\bigg\}\quad {\rm and}\cr
&u(x,t)\le\max_{Q_\varrho}u\le\mu_1-{\omega\over2^s}.\cr}
\leqno(5.10)
$$
Then (5.2) holds for any $\theta_1$. Let $\mu'_1>\mu_1-{\omega\over2^s}$. Then 
$H=\mu'_1-\left(\mu_1-{\omega\over2^r}\right)>\mu_1-{\omega\over2^s}-
\left(\mu_1-{\omega\over2^r}\right)={\omega\over2^r}-{\omega\over2^s}\ge
{\omega\over2^{s-1}}-{\omega\over2^s}={\omega\over2^s}$ for $r\le s-1$.

\noindent
Then Lemma 3.3 can be applied for cylinder $Q_\varrho$ and the level 
$k=\mu_1-{\omega\over2^r}$ if ${\omega\over2^r}\le\delta$. Hence Lemma 3.3 
implies
$$
{\rm meas}\bigg[B_\varrho\setminus 
A_{\mu_1-{\omega\over2^r}+{3\over4}H,\varrho}(t)\bigg]\ge{1\over18}
\varkappa_n\varrho^n
\leqno(5.11)
$$
for $t\in\left[t_0+{3\over2}\sigma\theta\varrho^2,t_0+
2\sigma\theta\varrho^2\right]$. Since
$$
\mu_1-{\omega\over2^r}+{3\over4}H\le\mu_1-{\omega\over2^r}+
{3\over4}{\omega\over2^r}=\mu_1-{\omega\over2^{r+2}},
$$
where we used that $H\le{\omega\over2^r}$, we obtain that
$$
{\rm meas}A_{\mu_1-{\omega\over2^r}+{3\over4}H,\varrho}(t)\ge{\rm meas}
A_{\mu_1-{\omega\over2^{r+2}}}(t).
$$
Hence, it follows from (5.11) the relation
$$
{\rm meas}[B_\varrho\setminus A_{\mu_1-{\omega\over2^{r+2}},\varrho}(t)]\ge
{1\over18}\varkappa_n\varrho^n,
\leqno(5.12)
$$
for $t\in\left[t_0+{3\over2}\sigma\theta\varrho^2,t_0+
2\sigma\theta\varrho^2\right]$.

\noindent
Since ${\omega\over2^r}\le\delta$ and $r\le s-1$, inequality (5.12) holds 
for $r\in\left[\log_2{\omega\over\delta},s-1\right]$.

\noindent
We have that $|\omega|\le2M$, where $M=\max_{Q_{\sigma\varrho}}u$ and we have 
also that $\log_2{2M\over\delta}\le{2M\over\delta}$. Therefore, we take 
$r\in\left[\left[{2M\over\delta}\right]+1,s-1\right]\subset
\left[\log_2{\omega\over\delta},s-1\right]$. 
We need the inequality (see [LSU, Ch. 2., (5.5)])
$$
(l-k){\rm meas}A_{l,\varrho}\le\beta
{\varrho^{n+1}\over{\rm meas}(B_\varrho\setminus A_{k,\varrho})}
\intop_{A_{k,\varrho}\setminus A_{l,\varrho}}|u_x|dx,
\leqno(5.13)
$$
where $\beta=\beta_1\varkappa_n^{1/n}$, 
$\beta_1={(1+n\varkappa_n)2^n\over n}$, 
$\varkappa_n=$ surface of $S^{n-1}$.

\noindent
Now we apply (5.13) to function $u(x,t)$ for the levels
$$
l=\mu_1-{\omega\over2^{p+1}},\quad k=\mu_1-{\omega\over2^p},\quad
t\in\bigg[t_0+{3\over2}\sigma\theta\varrho^2,t_0+
2\sigma\theta\varrho^2\bigg]
$$
and $p\in\left[\left[{2M\over\delta}\right]+3,s+1\right]$. 
Using (5.12) we obtain
$$
{\omega\over2^{p+1}}{\rm meas}A_{\mu_1-{\omega\over2^{p+1}},\varrho}(t)\le
18\beta\varkappa_n^{-1}\varrho\intop_{{\cal D}_p(t)}|u_x(x,t)|dx,
\leqno(5.14)
$$
where
$$
{\cal D}_p(t)=A_{\mu_1-{\omega\over2^p},\varrho}(t)\setminus
A_{\mu_1-{\omega\over2^{p+1}},\varrho}(t).
$$
Integrating (5.14) with respect to time in the time interval \\
$\left[t_0+{3\over2}\sigma\theta\varrho^2,t_0+2\sigma\theta\varrho^2\right]$ 
we obtain
$$\eqal{
&{\omega\over2^{p+1}}
\intop_{t_0+{3\over2}\sigma\theta\varrho^2}^{t_0+2\sigma\theta\varrho^2}
{\rm meas}A_{\mu_1-{\omega\over2^{p+1}},\varrho}(t)dt\cr
&\le18\beta\varkappa_n^{-1}\varrho
\intop_{t_0+{3\over2}\sigma\theta\varrho^2}^{t_0+2\sigma\theta\varrho^2}dt
\intop_{{\cal D}_p(t)}|u_x(x,t)|dx.\cr}
\leqno(5.15)
$$
Using that
$$
{\rm meas}Q_\varrho\bigg(\mu_1-{\omega\over2^{p+1}}\bigg)
\intop_{t_0+{3\over2}\sigma\theta\varrho^2}^{t_0+2\sigma\theta\varrho^2}
{\rm meas}A_{\mu_1-{\omega\over2^{p+1}},\varrho}(t)dt
$$
and taking square of (5.15) we derive
$$\eqal{
&\bigg({\omega\over2^{p+1}}\bigg)^2{\rm meas}^2Q_\varrho
\bigg(\mu_1-{\omega\over2^{p+1}}\bigg)\cr
&\le(18\beta\varkappa_n^{-1})^2\varrho^2\bigg|
\intop_{t_0+{3\over2}\sigma\theta\varrho^2}^{t_0+2\sigma\theta\varrho^2}dt
\intop_{{\cal D}_p(t)}|u_x(x,t)|dx\bigg|^2\cr
&\le(18\beta\varkappa_n^{-1})^2\varrho^2
\intop_{t_0+{3\over2}\sigma\theta\varrho^2}^{t_0+2\sigma\theta\varrho^2}dt
\intop_{{\cal D}_p(t)}u_x^2dx
\intop_{t_0+{3\over2}\sigma\theta\varrho^2}^{t_0+2\sigma\theta\varrho^2}
{\rm meas}{\cal D}_p(t)dt,\cr}
\leqno(5.16)
$$
where the Cauchy inequality was used. To estimate the first integral on the 
r.h.s. of (5.16) we use the inequality (3.10)
$$\eqal{
&\|\omega^{(k)}\|_{V_2^0(Q(\varrho-\sigma_1\varrho,\tau-\sigma_2\tau))}^2\le
\gamma\{[(\sigma_1\varrho)^{-2}+(\sigma_2\tau)^{-1}]
\|\omega^{(k)}\|_{L_2(Q(\varrho,\tau))}^2\cr
&\quad+c_1\delta^2\mu^{{2\over r}(1+\varkappa)}(k,\varrho,\tau)\},\quad 
{\rm where}\ \omega=\pm u(x,t).\cr}
\leqno(5.17)
$$
Replacing $Q(\varrho,\tau)$ by $Q_{\sigma\varrho}$ and 
$Q(\varrho-\sigma_1\varrho,\tau-\sigma_2\tau)$ by $Q_\varrho$ we obtain
$$\eqal{
&\intop_{t_0+{3\over2}\sigma\theta\varrho^2}^{t_0+2\sigma\theta\varrho^2}dt
\intop_{{\cal D}_p(t)}u_x^2(x,t)dx\le
\|u^{\left(\mu_1-{\omega\over2^p}\right)}\|_{V_2^0(Q_\varrho)}^2\cr
&\le\gamma\bigg\{\bigg[(\sigma\varrho)^{-2}+\bigg({3\over2}\sigma\theta
\varrho^2\bigg)^{-1}\bigg]
\|u^{\left(\mu_1-{\omega\over2^p}\right)}\|_{L_2(Q_{\sigma\varrho})}^2\cr
&\quad+c_1\delta^2\mu^{{2\over r}(1+\varkappa)}\bigg(\mu_1-{\omega\over2^p},
\sigma\varrho,2\sigma\theta\varrho^2\bigg)\bigg\}.\cr}
\leqno(5.18)
$$
To estimate the r.h.s. of (5.18) we examine
$$
\|u^{\left(\mu_1-{\omega\over2^p}\right)}\|_{L_2(Q_{\sigma\varrho})}^2=
\intop_{Q_{\sigma\varrho}}dxdt\bigg|u-\mu_1+{\omega\over2^p}\bigg|^2\le
\bigg({\omega\over2^p}\bigg)^2{\rm meas}Q_{\sigma\varrho},
$$
where meas 
$Q_{\sigma\varrho}=\varkappa_n(\sigma\varrho)^n2\sigma\theta\varrho^2$ 
and
$u-\mu_1\le0$ because $\mu_1=\max_{Q_{\sigma\varrho}}u$.

\noindent
Next we calculate
$$\eqal{
&\mu^{{2\over r}(1+\varkappa)}\bigg(\mu_1-{\omega\over2^p},\sigma\varrho,
2\sigma\theta\varrho^2\bigg)=\bigg[\intop_{t_0}^{2\sigma\theta\varrho^2+t_0}
{\rm meas}^{r\over q}A_{\mu_1-{\omega\over2^p},\sigma\varrho}(t)
dt\bigg]^{{2\over r}(1+\varkappa)}\cr
&\le\bigg[\intop_{t_0}^{2\sigma\theta\varrho^2+t_0}{\rm meas}^{r\over q}
B_{\sigma\varrho}dt\bigg]^{{2\over r}(1+\varkappa)}=
\{[\varkappa_n(\sigma\varrho)^n]^{r\over q}2\sigma\theta
\varrho^2\}^{{2\over r}(1+\varkappa)}.\cr}
$$
Using the above estimates in (5.18) yields
$$\eqal{
&\intop_{t_0+{3\over2}\sigma\theta\varrho^2}^{t_0+2\sigma\theta\varrho^2}dt
\intop_{{\cal D}_p(t)}u_x^2dx\cr
&\le\gamma\bigg\{\bigg[(\sigma\varrho)^{-2}+
\bigg({3\over2}\sigma\theta\varrho^2\bigg)^{-1}\bigg]
\bigg({\omega\over2^p}\bigg)^2\varkappa_n(\sigma\varrho)^n
2\sigma\theta\varrho^2\cr
&\quad+c_1\delta^2[[\varkappa_n(\sigma\varrho)^n]^{r\over q}2\sigma\theta
\varrho^2]^{{2\over r}(1+\varkappa)}\bigg\}.\cr}
\leqno(5.19)
$$
Since $\omega\ge2^s\varrho^{n\varkappa\over2}$ and 
$p\in\left[\left[{2M\over\delta}\right]+3,s+1\right]$ we obtain from (5.16) 
and (5.19) the inequality
$$\eqal{
&\bigg({\omega\over2^{p+1}}\bigg)^2{\rm meas}^2Q_\varrho
\bigg(\mu_1-{\omega\over2^{p+1}}\bigg)\cr
&\le(18\beta\varkappa_n^{-1})^2\varrho^2\gamma\bigg\{\bigg({1\over\sigma^2}+
{2\over3\sigma\theta}\bigg)\varrho^{-2}2\sigma\theta\varrho^2\varkappa_n
(\sigma\varrho)^n\bigg({\omega\over2^p}\bigg)^2\cr
&\quad+c_1\delta^2\varkappa_n^{{2\over q}(1+\varkappa)}
\theta^{{2\over r}(1+\varkappa)}\sigma^{\left(n{2\over q}+
{2\over r}\right)(1+\varkappa)}2^{{2\over r}(1+\varkappa)}
\varrho^{\left({2n\over q}+{4\over r}\right)(1+\varkappa)}\bigg\}\cr
&\quad\cdot
\intop_{t_0+{3\over2}\sigma\theta\varrho^2}^{t_0+2\sigma\theta\varrho^2}
{\rm meas}{\cal D}_p(t)dt.\cr}
\leqno(5.20)
$$
Using that 
$\varrho^{n\varkappa}\le\left({\omega\over2^s}\right)^2\le
\left({\omega\over2^p}\right)^2$ and
$\varrho^{\left({2n\over q}+{4\over r}\right)(1+\varkappa)}=
\varrho^n\varrho^{n\varkappa}$ we obtain from (5.20) the inequality
$$
\bigg({\omega\over2^{p+1}}\bigg)^2{\rm meas}^2Q_\varrho
\bigg(\mu_1-{\omega\over2^{p+1}}\bigg)\le c_2\varrho^{n+2}
\bigg({\omega\over2^p}\bigg)^2
\intop_{t_0+{3\over2}\sigma\theta\varrho^2}^{t_0+2\sigma\theta\varrho^2}
{\rm meas}{\cal D}_p(t)dt,
\leqno(5.21)
$$
where
$$\eqal{
c_2&=(18\beta\varkappa_n^{-1})^2\gamma\bigg[\bigg({1\over\sigma^2}+
{2\over3\sigma\theta}\bigg)2\sigma\theta\varkappa_n\sigma^n\cr
&\quad+4c_1\delta^2\varkappa_n^{{2\over q}(1+\varkappa)}
\theta^{{2\over r}(1+\varkappa)}
\sigma^{2\left({n\over q}+{1\over r}\right)(1+\varkappa)}\cdot
2^{{2\over r}(1+\varkappa)}\bigg].\cr}
$$
Hence, (5.21) implies
$$
{\rm meas}^2Q_\varrho\bigg(\mu_1-{\omega\over2^{p+1}}\bigg)\le4c_2\varrho^{n+2}
\intop_{t_0+{3\over2}\sigma\theta\varrho^2}^{t_0+2\sigma\theta\varrho^2}
{\rm meas}{\cal D}_p(t)dt.
\leqno(5.22)
$$
Since $p\in\left[\left[{2M\over\delta}\right]+3,s+1\right]$ then
$$
{\rm meas}Q_\varrho\bigg(\mu_1-{\omega\over2^{p+1}}\bigg)\ge{\rm meas}
Q_\varrho\bigg(\mu_1-{\omega\over2^s}\bigg)\quad {\rm for}\ \ p\le s-1.
\leqno(5.23)
$$
Employing (5.23) in (5.22) and summing with respect to $p$ from 
$\left[{2M\over\delta}\right]+3$ to $s-1$ we derive
$$\eqal{
&\bigg[s-\bigg[{2M\over\delta}\bigg]-4\bigg]{\rm meas}^2Q_\varrho
\bigg(\mu_1-{\omega\over2^s}\bigg)\cr
&\le4c_2\varrho^{n+2}
\intop_{t_0+{3\over2}\sigma\theta\varrho^2}^{t_0+2\sigma\theta\varrho^2}
{\rm meas}B_\varrho(t)dt\cr
&\le c_3\theta\varkappa_n\varrho^{2n+4}.\cr}
\leqno(5.24)
$$
Therefore for
$$
s=\bigg[{2M\over\delta}\bigg]+4+
\bigg[{c_3\theta\varkappa_n\over\theta_1^2}\bigg]
\leqno(5.25)
$$
we obtain (5.3). This concludes the proof.

\Remark{5.2.} 
In view of Remark 4.2 we know that $\theta_1$ is very small. Hence to minimize 
$s$ we assume that $\delta>2M$ and $\theta$ so small that 
${c_3\theta\varkappa_n\over\theta_1^2}<1$. Then
$$
s=4.
\leqno(5.26)
$$
Next we prove

\proclaim Lemma 5.3. (see Lemma 7.4 from [LSU, Ch. 2, Sect. 7]) 
For any $u\in{\cal B}_2(\Omega^T,M,\gamma,r,\delta,\varkappa)$ and for any 
cylinders $Q_{\varrho\over\sigma}$ and $Q_{\sigma\varrho}$, $\sigma>1$, 
with the same top and the axis the following inequalities hold:\\
either
$$
osc\{u,Q_{\varrho/\sigma}\}\le2^s\varrho^{n\varkappa\over2}
\leqno(5.27)
$$
or
$$
osc\{u,Q_{\varrho/\sigma}\}\le\bigg(1-{1\over2^s}\bigg)
osc\{u,Q_{\sigma\varrho}\},
\leqno(5.28)
$$
where $s=s(\theta_1)$ and $s(\theta_1)$ is described by (5.25) and (5.26).

\Proof 
Assume that
$$
\omega_1=osc\{u;Q_{\varrho/\sigma}\}>2^s\varrho^{n\varkappa\over2}.
$$
Then we also have that
$$
\omega=osc\{u,Q_{\sigma\varrho}\}>2^s\varrho^{n\varkappa\over2}.
$$
In view of Lemma 5.1 this implies inequalities either (5.2) or (5.3). 
Assume that (5.2) holds.
We apply Lemma 4.1 for function $u(x,t)$ defined in $Q_\varrho$ and for 
the level $k_0=\mu_1-{\omega\over2^{s-1}}$, where 
$\mu_1=\max_{Q_{\sigma\varrho}}u$ and either
$$
H=\max_{Q_\varrho}u-\bigg(\mu_1-{\omega\over2^{s-1}}\bigg)<
\varrho^{n\varkappa\over2}
\leqno(5.29)
$$
or
$$
{\rm meas}Q_{\varrho/\sigma}\bigg(k_0+{H\over2}\bigg)=0.
\leqno(5.30)
$$
Let us recall that the condition $k_0\ge\esssup_{Q_\varrho}u(x,t)-\delta$ 
from Lemma 4.1 hold. The first alternative (5.29) implies the inequality
$$
\max_{Q_{\varrho/\sigma}}u\le\mu_1-{\omega\over2^{s-1}}+\varrho^{\varkappa n\over2}
\le\mu_1-{\omega\over2^s}
\leqno(5.31)
$$
because ${\omega\over2^s}>\varrho^{n\varkappa\over2}$. Hence (5.31) yields
$$\eqal{
&\max_{Q_{\varrho/\sigma}}u-\min_{Q_{\varrho/\sigma}}u\le
\max_{Q_{\sigma\varrho}}u-\min_{Q_{\varrho/\sigma}}u-{\omega\over2^s}\cr
&\le\max_{Q_{\sigma\varrho}}u-\min_{Q_{\sigma\varrho}}u-{\omega\over2^s},\cr}
\leqno(5.32)
$$
where we used that $\min_{Q_{\varrho/\sigma}}u\ge\min_{Q_{\sigma\varrho}}u$.
Therefore (5.32) takes the form
$$
osc\{u,Q_{\varrho/\sigma}\}\le\bigg(1-{1\over2^s}\bigg)osc
\{u,Q_{\sigma\varrho}\},
\leqno(5.33)
$$
so (5.28) holds. Let us consider the second alternative (5.30). It implies
$$\eqal{
&\max_{Q_{\varrho/\sigma}}u\le\mu_1-{\omega\over2^{s-1}}+{H\over2}=
\mu_1-{\omega\over2^{s-1}}+{1\over2}\bigg[\max_{Q_\varrho}u-
\bigg(\mu_1-{\omega\over2^{s-1}}\bigg)\bigg]\cr
&=\mu_1-{\omega\over2^{s-1}}+{1\over2}(\max_{Q_\varrho}u-
\max_{Q_{\sigma\varrho}}u)+{\omega\over2^s}\le\mu_1-{\omega\over2^s},\cr}
\leqno(5.34)
$$
because $\max_{Q_\varrho}u\le\max_{Q_{\sigma\varrho}}u$.

\noindent
Similarly as above (5.34) implies (5.28). This proves the lemma.

\noindent
From Lemmas 5.3 and 2.5 we have

\proclaim Theorem 5.4. 
Let $u(x,t)\in{\cal B}_2(\Omega^T,M,\gamma,r,\delta,\varkappa)$ and 
$Q_{\varrho_0}=Q(\varrho_0,\theta\varrho_0^2)\subset\Omega^T$. Then for any 
cylinder $Q_\varrho=Q(\varrho,\theta\varrho^2)$ with the same top and the 
axis as $Q_{\varrho_0}$ we obtain
$$
osc\{u,Q_\varrho\}\le c\varrho^\alpha\varrho_0^{-\alpha}
\leqno(5.35)
$$
where
$$
\alpha=\min\bigg\{-\ln_{\sigma^2}\bigg(1-{1\over2^s}\bigg),
{n\varkappa\over2}\bigg\},\quad 
c=\sigma^{2\alpha}\max\bigg\{2M,2^s\varrho_0^{n\varkappa\over2}\bigg\}.
$$

\Remark{5.5.} 
In our case $s=4$, ${n\varkappa\over2}={1\over2}$. Therefore
$$
\alpha=\min\left\{-\ln_{\sigma^2}{15\over16},{1\over2}\right\}.
$$
Hence, we can choose $\sigma$ such that $\ln_{\sigma^2}{16\over15}={1\over2}$ 
so $\sigma={16\over15}$.

\noindent
Hence Theorem 5.4 and Remark 5.5 imply the Main Theorem.

\section{References}

\item{LL.} Landau, L.; Lifshitz, E.: Hydrodynamics, Nauka, Moscow 1986 
(in Russian).

\item{LSU.} Ladyzhenskaya, O. A.; Solonnikov, V. A.; Uraltseva, N. N.: 
Linear and quasilinear equations of parabolic type, Nauka, Moscow 1967 
(in Russian).

\item{Z1.} Zaj\c aczkowski, W. M.: Global regular axially symmetric solutions 
to the Navier-Stokes equations in a periodic cylinder

\item{Z2.} Zaj\c aczkowski, W. M.: Global special regular solutions to the 
Navier\--Stokes equations in a cylindrical domain under boundary slip 
condition, Gakuto Intern. Ser., Math. Sc. Appl. 21 (2004), pp. 188.

\item{Z3.} Zaj\c aczkowski, W. M.: A priori estimate for axially symmetric 
solutions to the Navier-Stokes equations near the axis of symmetry.

\bye